\documentclass[10pt]{article}
\usepackage{amsmath,amsthm,amsfonts,amscd,amssymb} 
\usepackage{verbatim}   
\usepackage{eucal} 
\usepackage{longtable} 
\usepackage{makeidx}
\usepackage[T1]{fontenc}
\usepackage[utf8]{inputenc}
\usepackage{authblk}

\usepackage{hyperref}

\newcommand{\field}[1]{\mathbb{#1}}

\newcommand{\R}{\field{R}}
\newcommand{\C}{\field{C}}

\newcommand{\be}{\begin{equation}}
\newcommand{\bel}[1]{\begin{equation}\label{#1}}
\newcommand{\ee}{\end{equation}}

\newcommand{\ba}{\begin{eqnarray}}
\newcommand{\ea}{\end{eqnarray}}

\newcommand{\bas}{\begin{eqnarray*}}
\newcommand{\eas}{\end{eqnarray*}}

\newcommand{\thm}{\begin{theorem}}
\newcommand{\thml}[1]{\begin{theorem}\label{#1}}
\newcommand{\mht}{\end{theorem}}
\newcommand{\cnj}{\begin{conjecture}}
\newcommand{\cnjl}[1]{\begin{conjecture}\label{#1}}
\newcommand{\jnc}{\end{conjecture}}
\newcommand{\dfn}{\begin{definition}}
\newcommand{\dfnl}[1]{\begin{definition}\label{#1}}
\newcommand{\nfd}{\end{definition}}
\newcommand{\pro}{\begin{proposition}}
\newcommand{\prol}[1]{\begin{proposition}\label{#1}}
\newcommand{\orp}{\end{proposition}}
\newcommand{\crl}{\begin{corollary}}
\newcommand{\crll}[1]{\begin{corollary}\label{#1}}
\newcommand{\lrc}{\end{corollary}}
\newcommand{\lmm}{\begin{lemma}}
\newcommand{\lmml}[1]{\begin{lemma}\label{#1}}
\newcommand{\mml}{\end{lemma}}
\newcommand{\prf}{\begin{proof}}
\newcommand{\prfl}[1]{\begin{proof}\label{#1}}
\newcommand{\frp}{\end{proof}}
\newcommand{\axi}{\begin{axiom}}
\newcommand{\axil}[1]{\begin{axiom}\label{#1}}
\newcommand{\ixa}{\end{axiom}}
\newcommand{\rmk}{\begin{remark}}
\newcommand{\rmkl}[1]{\begin{remark}\label{#1}}
\newcommand{\kmr}{\end{remark}}
\newcommand{\exa}{\begin{example}}
\newcommand{\exal}[1]{\begin{example}\label{#1}}
\newcommand{\axe}{\end{example}}
\newcommand{\alg}{\begin{algorithm}}
\newcommand{\gla}{\end{algorithm}}
\newcommand{\nte}{\begin{note}}
\newcommand{\ntel}[1]{\begin{note}\label{#1}}
\newcommand{\etn}{\end{note}}
\newcommand{\app}{\begin{application}}
\newcommand{\appl}[1]{\begin{application}\label{#1}}
\newcommand{\ppa}{\end{application}}
\newcommand{\mat}{\begin{matrix}}
\newcommand{\tam}{\end{matrix}}
\newcommand{\smm}{\begin{summary}}
\newcommand{\mms}{\end{summary}}

\newcommand{\teo}{\begin{teorema}}
\newcommand{\teol}[1]{\begin{teorema}\label{#1}}
\newcommand{\oet}{\end{teorema}}
\newcommand{\df}{\begin{definicion}}
\newcommand{\dfl}[1]{\begin{definicion}\label{#1}}
\newcommand{\fd}{\end{definicion}}
\newcommand{\por}{\begin{proposicion}}
\newcommand{\porl}[1]{\begin{proposicion}\label{#1}}
\newcommand{\rop}{\end{proposicion}}
\newcommand{\cor}{\begin{corolario}}
\newcommand{\corl}[1]{\begin{corolario}\label{#1}}
\newcommand{\roc}{\end{corolario}}
\newcommand{\lem}{\begin{lema}}
\newcommand{\leml}[1]{\begin{lema}\label{#1}}
\newcommand{\mel}{\end{lema}}
\newcommand{\pru}{\begin{prueba}}
\newcommand{\prul}[1]{\begin{prueba}\label{#1}}
\newcommand{\urp}{\end{prueba}}
\newcommand{\axa}{\begin{axioma}}
\newcommand{\axal}[1]{\begin{axioma}\label{#1}}
\newcommand{\xa}{\end{axioma}}
\newcommand{\nta}{\begin{nota}}
\newcommand{\ntal}[1]{\begin{nota}\label{#1}}
\newcommand{\atn}{\end{nota}}
\newcommand{\eje}{\begin{ejemplo}}
\newcommand{\ejel}[1]{\begin{ejemplo}\label{#1}}
\newcommand{\je}{\end{ejemplo}}
\newcommand{\api}{\begin{aplicacion}}
\newcommand{\apil}[1]{\begin{aplicacion}\label{#1}}
\newcommand{\ipa}{\end{aplicacion}}

\newcommand{\ali}{\begin{align}}
\newcommand{\ila}{\end{align}}
\newcommand{\enu}{\begin{enumerate}}
\newcommand{\une}{\end{enumerate}}
\newcommand{\arr}{\begin{array}}
\newcommand{\rra}{\end{array}}
\newcommand{\eqa}{\begin{eqnarray}}
\newcommand{\aqe}{\end{eqnarray}}
\newcommand{\equ}{\begin{equation}}
\newcommand{\uqe}{\end{equation}}
\newcommand{\subq}{\begin{subequations}}
\newcommand{\qbus}{\end{subequations}}

\newtheorem{theorem}{Theorem}[section]
\newtheorem{lemma}[theorem]{Lemma}
\newtheorem{proposition}[theorem]{Proposition}
\newtheorem{corollary}[theorem]{Corollary}
\newtheorem{conjecture}[theorem]{Conjecture}

\theoremstyle{definition}
\newtheorem{definition}{Definition}[section]
\newtheorem{axiom}[definition]{Axiom}

\theoremstyle{remark}
\newtheorem{remark}{Remark}[section]
\newtheorem{example}[remark]{Example}
\newtheorem{application}[remark]{Application}
\newtheorem{algorithm}[remark]{Algorithm}
\newtheorem{note}[remark]{Note}
\newtheorem{summary}[remark]{Summary}

\newtheorem{teorema}[theorem]{Teorema}
\newtheorem{lema}[theorem]{Lema}
\newtheorem{proposicion}[theorem]{Proposici\'on}
\newtheorem{corolario}[theorem]{Corolario}

\newtheorem{definicion}[definition]{Definici\'on}
\newtheorem{axioma}[definition]{Axioma}

\newtheorem{nota}[remark]{Nota}
\newtheorem{ejemplo}[remark]{Ejemplo}
\newtheorem{aplicacion}[remark]{Aplicaci\'on}

\begin{document}

\title{Multiplicative operators in the spectral problem of  integrable systems}

\author[$\dagger$]{ A. ESP\'INOLA-ROCHA}

\affil[$\dagger$]{{\it Departamento de Ciencias B\'asicas. 

Universidad Aut\'onoma Metropolitana-Azcapotzalco.

Av. San Pablo No. 180. 
Col. Reynosa-Tamaulipas. 
C.P. 02200.
Azcapotzalco, 
CDMX,
M\'exico }

{\tt Email:} {jaer.azc.uam.mx@gmail.com}}

\author[$\ddagger$]{F. X. PORTILLO-BOBADILLA}

\affil[$\ddagger$]{{\it Universidad Aut\'onoma de la Ciudad de M\'exico.

Calzada Ermita Iztapalapa No. 4163, Lomas de Zaragoza, C.P. 09620
 Iztapalapa, 
CDMX,  M\'exico}

{\tt Email:} {francisco.portillo@uacm.edu.mx}}


\maketitle


\begin{abstract}
We consider the spectral problem of the Lax pair associated to  periodic integrable partial differential equations. We assume this spectral problem to be a polynomial of degree 
$d$ in the spectral parameter $\lambda$. From this assumption, we find the conservation laws as well as the hyperelliptic curve required to solve the periodic inverse problem. A recursion formula is developed, as well 
as $d$ additional conditions which give additional information to integrate the equations under consideration.
We also include two examples to show how the techniques developed work. For the Korteweg-deVries  (KdV) equation, the degree of the multiplicative equation is $d=1$. Hence, we only have one condition and one recursion formula. The condition gives in each degree of the recursion the conserved  densities for KdV equation, recovering the Lax hierarchy. For the Nonlinear Schr\"odinger (NLS) equation, the degree of the multiplicative operator is $d=2$. Hence, we have a couple of conditions that we use to deduce conserved density constants and the Lax hierarchy for such equation. Additionally, we explicitely write down the hyper-elliptic curve associted to the NLS equation. Our approach can be use for other completely  integrable differential equations, as long as we have a polynomial multiplicative operator associated to them.

 \end{abstract}

\section*{Introduction}

John Scott Russell observed for the first time a {\em "soliton"} (although he named it a {\em "great wave of translation"})  in the 1830s in a channel in Edinburgh, Scotland. Russell performed several experiments and  reported his results in 1834 (see reference in \cite{scott}). But it was not until 1895 that a partial differential equation (PDE) was developed by two dutch mathematicians, D.J. Korteweg and G. de Vries  \cite{KdV_seminal}, to describe Russell's {\em "great wave"}. In their paper, the authors found the one-soliton solution, a sech$^2$-profile describing the wave observed by Russell, as well as a periodic solution of cnoidal-type. This equation is now known as the Korteweg-deVries (KdV) equation, 
\begin{equation}\label{kdv_equation}
\frac{\partial q}{\partial t} + 6q \frac{\partial q}{\partial x}  +   \frac{\partial^3 q}{\partial x^3}    = 0,
\end{equation}
where $x$ denotes de spacial variable along the channel, $t$ the temporal variable, and $q= q(x,t)$ the elevation of the wave at the position $x$ and at time $t$.

In 1967, C. Gardner, J. Greene, M. Kruskal and R. Miura \cite{ggkm} found a method to exactly solve the KdV equation, the {\em inverse scattering method} (or {\em transform}) (IST). A year later, P. Lax generalized this idea stating that all the equations that possess  a pair of {\em linear} operators to express the equations can be solved by the IST method. Since then, the pair of operators have been  called the {\em Lax pair} \cite{lax_1}. In 1972, following Lax \cite{lax_1},  V.E. Zakharov and A.B. Shabat were able to integrate the nonlinear Schr\"odinger equation (NLS)  \cite{zajarov}, 
\begin{equation}\label{nls_equation}
i  \frac{\partial q}{\partial t} =  -\frac{1}{2}  \frac{\partial^2 q}{\partial x^2}    -  \sigma|q|^2q,
\end{equation}
where $|q|^2$ represents the intensity of an electromagnetic field inside an optical fiber and $\sigma = \pm1$ represents the focusing or defocusing nature  of the fiber.
In 1974, another equation, the sine-Gordon equation, was completely integrated  \cite{akns} and the method could be further generalized \cite{akns2}.  After that, many other equations, continuous and discrete, were found to be  completely integrable (see \cite{scott, drazin, nmpz, ablowitz, ablowitz2, newell} and references therein). 

The IST was originally developed  in the real line but earlier in the 1970s, P. Lax worked  the KdV equation in a periodic domain \cite{lax_2}. Independently, and about the same years, several authors  such as B.A. Dubrovin, A.P. Its,  I.M. Krichever, V.B. Matveev and S.P. Novikov were also working  the periodic KdV equation using the  techniques of algebraic geometry \cite{dubrovin_1, dubrovin_2, dubrovin_novikov_1, dubrovin_novikov_2, dubrovin_matveev_novikov, krichever_1, krichever_2, krichever_3, krichever_novikov, novikov, its_matveev}. A good survey of the topic can be found in \cite{dubrovin_1, krichever_1}. In the present paper, we  follow the ideas in \cite{dubrovin_1} and \cite{flaschka_1}.

At this point,  almost all  the Lax pairs found for  new integrable systems were matrix operators (as an example, the Lax pair for the KdV equation is scalar, although it  also could be written as pair of $2\times2$ matrix operators). However, in the work of Kamchatnov, Kraenkel and Umarov  \cite{Kamchatnov_Kraenkel}  the authors found  the corresponding scalar Lax pairs to the usual $2\times2$ Lax matrix representations. In particular, applications can be found in  \cite{Kamchatnov_Kraenkel_Umarov_1,  Kamchatnov_Kraenkel_Umarov_2, Kamchatnov_1, Kamchatnov_2}. The authors applied  this technique \cite{Kamchatnov_Kraenkel_Umarov_1, Kamchatnov_Kraenkel_Umarov_2} to the nonlinear Schr\"odinger equation, whose associated spectral problem is a polynomial of degree 2 in the spectral parameter, to find trains of solitons using the Bohr-Sommerfeld quantization rule, just as it is done in the quantum case for the Schr\"odinger equation \cite{landau_lifshitz_3}.

In this paper, we consider the  spectral problem associated to certain nonlinear evolution equation, {\em i.e.}, the problem associated to the Lax operator of the spatial evolution. We assume that this   problem can be written in  scalar form  and that the corresponding Lax operator is a {\em multiplicative operator}.
Even more, we assume that this operator is a polynomial  of degree $d$ in the spectral parameter, $d$ being a positive integer. We consider the work of Dubrovin {\em et.al.} \cite{dubrovin_1, dubrovin_matveev_novikov}, but follow Flaschka \cite{flaschka_1}. Under this approach, we find the corresponding general hyperelliptic curve, a recursion formula to find the coefficients of the $n$-soliton solution,  the conserved densities  of motion, as well as the equations in the Lax hierarchy. The integrals of motion appear as coefficients in the polynomial expansion of the squared eigenfunction solution  and  the Lax hierarchy  appears in extra conditions of solvability for the third order differential equation of the squared eigenfunction. As particular examples, we consider the KdV  (\ref{kdv_equation}) ($d= 1$) and the NLS  (\ref{nls_equation})  ($d=2$) equations. As far as the authors know, there is no approach similar to the methodology presented here for $d$ being an arbitrary positive integer. 

To see how this works, consider a solution $y(x,t)$ of the Schr\"odinger-type equation,
\equ\label{inicial}
\frac{\partial^2y}{\partial x^2}=\hat{\mathcal{L}}y,
\uqe
where $\hat{\mathcal{L}}$ is a multiplicative linear operator, polynomial of degree $d$ in the spectral parameter $\lambda$. For the KdV equation (\ref{kdv_equation}), $\hat{\mathcal{L}} = \lambda -q$, {\em i.e.}, $d=1$.
Since $\hat{\mathcal{L}}$ is a multiplicative operator, multiplying (\ref{inicial}) by $y$ 
results into
\equ\label{inicial2}
y\frac{\partial^2y}{\partial x^2}=\hat{\mathcal{L}}y^2
\uqe
Now, setting  $\phi=y^2$ (this is  the {\em squared eigenfunction} of $\hat{\mathcal{L}}$), taking the derivative with respect to $x$, $\phi'=2yy'$ and squaring, we obtain:
$
(y')^2=\left(\phi'\right)^2/ (4\phi).            
$
Similarly, taking a second derivative of $\phi$, we have
$
\phi'' 
=2(y')^2+2yy''
= \left(\phi'\right)^2/(2\phi)+2yy''.                   
$
Hence,
$
  yy''=\frac{1}{2}\phi''- (\phi')^2/(4\phi),           
$
and equation (\ref{inicial2}) becomes
\equ\label{eq00}
\frac{1}{2}\phi''-\frac{(\phi')^2}{4\phi}=\hat{\mathcal{L}}\phi.
\uqe
Taking derivatives on both sides, we have
\equ\label{rec:1}
\frac{1}{2}\phi'''-\left(\frac{(\phi')^2}{4\phi}\right)_x=\left[\hat{\mathcal{L}}\phi\right]_x
\uqe

A 
 straight-forward computation, using equation (\ref{eq00}), shows that
$
\left((\phi')^2/(4\phi)\right)_x 
=   (\phi'/\phi)\left(\frac{1}{2}\phi''-  (\phi')^2/(4\phi)\right)
=  (\phi'/\phi) \hat{\mathcal{L}}\phi=\hat{\mathcal{L}}\phi' ,
$
since $\hat{\mathcal{L}}$ is a multiplicative operator. 

Thus, equation (\ref{rec:1}) turns to be
%
 the following third order differential equation for the squared eigenfunction $\phi$:
\equ\label{rec:4}
\phi'''-2\hat{\mathcal{L}}_x\phi-4\hat{\mathcal{L}}\phi'=0.
\uqe

To keep computations compact, we define the following bilinear operator:
\equ\label{bilinear:operator}
\langle \psi, \phi \rangle 
:= (\psi\phi)_x+\psi\phi_x
= \psi_x\phi+2\psi\phi_x,
\uqe
so that equation (\ref{rec:4}) can be written as:
\equ\label{recursion:1}
\phi'''-2\langle\hat{\mathcal{L}} , \phi \rangle=0.
\uqe
Equation (\ref{recursion:1}) is the differential equation to be solved for the squared eigenfunction $\phi$. {\em Notice
that this equation holds no matter what  the degree of the operator $\hat{\mathcal{L}}$ is}. Here it is where all of our argument is sustained. 

Then, the paper is organized as follows. In section \ref{sect:recursion}, it is found a recursion formula to determine the coefficients of the $n$-degree polynomial $\phi_n(x;\lambda)$, the squared eigenfunction, called here $n$-soliton solutions. Additionally, $d$ extra conditions of solvability are found. In section \ref{sect:hyperelliptic_curve}, we explicitly construct the hyperelliptic curve associated to the spectral problem, equation (\ref{inicial}). Here, a polynomial of degree $2n+d$ is obtained. All the coefficients of the polynomial are computed in terms of the coefficients of the squared eigenfunction  $\phi_n(x;\lambda)$ and the operator $\hat{\mathcal{L}}$ (considered as polynomials in $\lambda$); therefore, they implicitly  depend on the solution $q(x)$. Also, a set of  differentials on the hyperelliptic curve is constructed. This set of differentials is essential for the integration and the inverse problem associated to these type of equation, but it will not be done here. The set of differentials is also found in  \cite{newell, espinola_portillo, miller}. Section \ref{sect:examples} contains the KdV and NLS equations as examples of the theory developed in previous sections. A conclusion section comes at the end of the paper.





\subsubsection*{Remark.} It is important to mention here that, as in \cite{lax_1, novikov}, we also work the stationary case of integrable equations only. Therefore, we just consider the spectral problem and the corresponding element of the Lax pair, not the time evolution  associated  to the second element of the Lax pair. 


\subsubsection*{Notation.} In this work, $f'(x)$ and $f_x(x)$ represent the derivative of $f(x)$ with respect to  $x$ (even if $f$ depends on more variables). For the KdV equation (subsection \ref{subsect:kdv}), we prefer $q'$, $q''$,$\ldots$; for the derivatives of $q$ of successive order with respect to $x$. But, in the NLS equation (subsection \ref{subsect:nls}), we use $q_x$, $q_{xx}$,$\ldots$; for the derivatives of $q$. In this last case, the derivatives of $E$ and $F$ (see equation \ref{nls:equation:E:F}), are better denoted $E'$, $E''$, $F$, $F''$, etc.

\subsubsection*{Nomenclature.} A {\em squared eigenfunction, polynomial of degree $n$ in the variable $\lambda$}, $\phi_n(x)$,  will be  called in this paper an \textbf{$n$-soliton solution}.  See definition \ref{normalized_n_soliton}. (It is usual to call {\em "solitons"} to exact solutions of integrable partial differential equations, although we will not name them this way here).


\section{Recursion formul\ae \hspace{0.08 cm} and conditions of solvability}\label{sec:recursion_formulas_conditions}\label{sect:recursion}

Assume 
$\hat{\mathcal{L}}$ is a multiplicative operator,  polynomial  of degree $d$ in $\lambda$, with coefficients depending on $x$,
\be \label{L_degree_d}
\hat{\mathcal{L}}(x;\lambda) = \mathcal{L}_0 \lambda^d  +  \mathcal{L}_1(x) \lambda^{d-1} + \dots + \mathcal{L}_d(x) = \sum\limits_{j = 0}^d \mathcal{L}_j(x) \lambda^{d-j},
\ee
and $\phi_n$ is an squared eigenfunction,  and {\em  $n-\text{soliton}$} solution,  (associated to $\hat{\mathcal{L}}(x;\lambda)$), also polynomial  in $\lambda$, but of degree $n$,
\be \label{squared_eigenfuntion_n}
\phi_n(x;\lambda) = A_0 \lambda^n  +  A_1(x) \lambda^{n-1} + \dots + A_n(x) = \sum\limits_{i = 0}^n A_i(x) \lambda^{n-i},
\ee
where $A_0 \neq 0 $ and $\mathcal{L}_0 \neq 0$ are $x$-independent, {\em i.e.} are constants (that might depend on $t$ in the non-stationary case), $A_i$ and $\mathcal{L}_j$  are functions of $x$ only (and they might also depend on $t$ in the non-stationary case). And $\lambda$ is a constant (in principle, it could be  a function of $t$ only in the non-stationary case).

Recursion {formul\ae} for the squared eigenfunctions $\phi_n$ and conditions for the solutions are obtained when solving the equation
\equ\label{recursion:2}
\phi_n''' - 2\langle \hat{ \mathcal{L}}, \phi_n  \rangle = 0,
\uqe
by setting to zero each coefficient of the polynomial. Notice that the term $2\langle \hat{ \mathcal{L}}, \phi_n  \rangle $ turns to be a polynomial of at most degree $(n+d)$. From here, it is where the recursion formul\ae \; and  the $d$ extra conditions arise. 

In fact, by bilinearity of the product $\langle  \ , \ \rangle$, equation (\ref{bilinear:operator}) is 
\equ\label{bilinearity_sum}
\begin{split}
\phi_n''' - 2 \langle \hat{ \mathcal{L}}, \phi_n  \rangle  = &  \sum\limits_{i=0}^n A_i'''\lambda^{n-i} - 2\left\langle   \sum\limits_{i=0}^d \mathcal{L}_i\lambda^{d-i}  ,   \sum\limits_{j=0}^n A_j\lambda^{n-j} \right\rangle \\
= &   \sum_{i=0}^n A_i'''\lambda^{n-i}  - 2 \sum_{k= 0}^{d+n}\left(  \sum_{\scriptsize\begin{matrix} i+j = k\\ i,j\geq 0\end{matrix}}  \langle  \mathcal{L}_i,   A_j  \rangle \right)\lambda^{d+n-k}
\end{split}
\uqe

If we set all coefficients of (\ref{bilinearity_sum}) equal to zero, we find relations among the coefficients $A_i$ of the solutions $\phi_n(x;\lambda)$.

\underline{For $k = 0$}, the fact that $\mathcal{L}_0$ and $A_0$ are constants implies that $\langle    \mathcal{L}_0 , A_0   \rangle = 0$ and, hence, there is no term of degree $(d+n)$. 

\underline{For $k = 1$}, it follows
$
\langle  \mathcal{L}_0 , A_1 \rangle   + \langle \mathcal{L}_1 , A_0 \rangle  = 0,
$
and solving for $A_1$,
$$
A_1 = -\;\frac{1}{2\mathcal{L}_0}\int   \langle  \mathcal{L}_1 , A_0 \rangle   \; dx.
$$

\underline{For $k = 2$},  we can check that  $\langle  \mathcal{L}_0 , A_{2} \rangle   +     \langle  \mathcal{L}_1 , A_{1} \rangle +     \langle  \mathcal{L}_2 , A_{0} \rangle   =    
\langle  \mathcal{L}_0 , A_{2} \rangle   + \sum\limits_{i = 1}^{2}    \langle  \mathcal{L}_i , A_{2-i} \rangle   =   0$.

In general, \underline{for $k \leq d$}  proceeding recursively, we obtain $\langle  \mathcal{L}_0 , A_{k} \rangle   + \sum\limits_{i = 1}^{k}    \langle  \mathcal{L}_i , A_{k-i} \rangle   =   0$. 
Solving for $A_{k},$
$$
A_{k}  =    -\;\frac{1}{2\mathcal{L}_0}\int      \sum\limits_{i = 1}^{k}    \langle  \mathcal{L}_i , A_{k-i} \rangle       \; dx.
$$


Starting at \underline{$k = d+1$}, the polynomial 
$
\phi_n'''=\sum\limits_{i=1}^{n}A_i''' \lambda^{n-i}
$
of degree $n-1$ (remember $A_0=$ constant) should be now taken into account. 
For example, for \underline{$k = d+1$}, we obtain the formula:
$$
A_{d+1}  =  -\;\frac{1}{2\mathcal{L}_0}\left( \int      \sum\limits_{i = 1}^{d}    \langle  \mathcal{L}_i , A_{d+1-i} \rangle       \; dx  - \frac{1}{2}A_1''\right).
$$

Hence, the \underline{general recursion formula} to compute the coefficients of $\phi_n(x;\lambda)$ is:
\be\label{recursion_formula}
A_{k}  =  -\;\frac{1}{2\mathcal{L}_0}\left( \int      \sum\limits_{i = 1}^{d}    \langle  \mathcal{L}_i , A_{k-i} \rangle       \; dx  - \frac{1}{2}A_{k-d}''\right),
\ee
which holds for $k = 1, 2, \dots, n$ (considering $A_{-l} = 0$ for $l>0$).

Now, in order to have a solution to equation (\ref{recursion:2}), a few extra conditions are needed. They come from imposing the coefficients of $\lambda^{d-1},\lambda^{d-2}, \dots, \lambda, 1,$ in (\ref{bilinearity_sum}) to be zero. 

Setting the coefficient of $\lambda^{d-1}$ equal to zero,  we obtain the first condition:
$$
A_{n-d+1}''' - 2  \sum\limits_{i = 1}^{d}    \langle  \mathcal{L}_i , A_{n+1-i} \rangle   = 0.
$$
 
In general, the  \underline{$d$ conditions} to solve the system are
\be\label{d_conditions}
\mathcal{A}_{n,s}:=A_{n-s}''' - 2  \sum\limits_{i = d-s}^{d}    \langle  \mathcal{L}_i , A_{d + n-s-i} \rangle   = 0,
\ee
for $0\leq s \leq d-1$, where $A_{n-s} = 0$ for values of $s$ such that  $n-s < 0$. (This includes the cases $n\leq d$ and $d<n$).

There is an alternative form to express this condition by the use of  the recursion formula (\ref{recursion_formula}). In fact, writing:
\begin{equation*}
\begin{split}
\mathcal{A}_{n,s}&=A_{n-s}''' - 2  \sum\limits_{i = d-s}^{d}    \langle  \mathcal{L}_i , A_{d + n-s-i} \rangle  \\ &= A_{n-s}''' - 2  \sum\limits_{i = 1}^{d}\langle  \mathcal{L}_i , A_{d + n-s-i} \rangle +2\sum\limits_{i = 1}^{d-s-1} \langle  \mathcal{L}_i , A_{d + n-s-i} \rangle  \\
&= 4\mathcal{L}_0\left\lbrace -\frac{1}{2\mathcal{L}_0}\left(\sum\limits_{i = 1}^{d}\langle  \mathcal{L}_i , A_{d + n-s-i} \rangle-\frac{1}{2}A_{n-s}'''\right)\right\rbrace+2\sum\limits_{i = 1}^{d-s-1} \langle  \mathcal{L}_i , A_{d + n-s-i} \rangle \\
&= 4\mathcal{L}_0A'_{n+d-s}+2\sum\limits_{i = 1}^{d-s-1} \langle  \mathcal{L}_i , A_{d + n-s-i} \rangle ,  \quad \text{by the recursion (\ref{recursion_formula})},
\end{split}
\end{equation*}
therefore obtaining  the \underline{alternative form for the conditions}: 
\begin{equation}\label{d_conditions_alternative}
\mathcal{A}_{n,s} = 4\mathcal{L}_0A'_{n+d-s}+2\sum\limits_{i = 1}^{d-s-1} \langle  \mathcal{L}_i , A_{d + n-s-i} \rangle = 0 .
\end{equation}
We will  use this alternative form to express Theorems \ref{all_KdV_N_solitons}, \ref{NLS_theorem} and \ref{all_NLS_N_solitons} in the examples for the KdV and the NLS equations. 

We must point out that we can construct all the solutions to equation (\ref{recursion:2}) using the general recursion formula subject to the $d$ conditions.
\dfn \label{normalized_n_soliton}
Define the \emph{\bf normalized $n$-soliton solution 
associated to the linear operator $\hat{\mathcal{L}}$} to be the solution to (\ref{recursion:1}) obtained by the general recursion formula, assuming that $A_0=-2\mathcal{L}_0$, and that all the constants of integration in (\ref{recursion_formula}) are zero. Also, name $\psi_n$ a polynomial  solution  to (\ref{recursion:1}) of degree $n$ in $\lambda$, an \emph{\bf $n$-soliton solution associated to the linear operator $\hat{\mathcal{L}}$}. 
\nfd
Thus, the following theorem holds.




\thm\label{all_N_solitons}
\enu
\item[A)] Each $n$-soliton $\psi_n = \sum\limits_{k=0}^n B_k(x) \lambda^{n-k}$ can be written as a linear combination of the normalized $n$-solitons: $\phi_n$, $\phi_{n-1}$, $\ldots$, $\phi_0$.
\item[B)] The linear combination
$$\psi_n=K_0\phi_n+K_{1}\phi_{n-1}+\ldots+K_n\phi_0$$ (with $K_i$ constant for $0\leq i\leq n$ and $K_0\neq 0$) is an $n$-soliton solution if, and only if, it satisfies the following conditions
\equ\label{linear_combination_of_conditions}
\mathcal{B}_{n,s}=\sum_{l=0}^n K_l\mathcal{A}_{n-l,s}=0,  \quad   \text{ for }  \quad 0\leq s \leq d-1 ,
\uqe
where $\mathcal{A}_{n-l,s}$ and $\mathcal{B}_{n,s}$ are given as in equation (\ref{d_conditions})  for $\phi_n$ and $\psi_n$, respectively. Since $\mathcal{A}_{n-l,s} = 0$ is the {$s^{\text{th}}$ condition} for the normalized soliton  $\phi_{n-l}$, then $\mathcal{B}_{n,s} = 0$ is the {$s^{\text{th}}$ condition} for the  soliton  $\psi_{n}$, {\em i.e.}, similar conditions hold for $\psi_n$.
\une
\mht

\prf
\enu

\item[A)]

Let $A_k$ be the coefficient of the normalized soliton corresponding to the  monomial $\lambda^{n-k}$ for $\phi_n$, with $A_0 \neq 0$. 
Let 
$$\psi_n(x;\lambda)  =\sum_{k=0}^n B_k(x) \lambda^{n-k}$$
be an $n$-soliton of degree $n$, with $B_0 \neq 0$. .

Since $A_0\neq 0$ and $B_0\neq 0$ are constants, then  $K_0:=B_0A^{-1}_0$ is a constant. Hence, using 
the recursion formula on the coefficients of $\psi_n$, we compute:
\begin{equation*}
\begin{split}
B_1 &=-\frac{1}{2\mathcal{L}_0}\int \langle \mathcal{L}_1, B_0\rangle \; dx\\
&=K_0\left( -\frac{1}{2\mathcal{L}_0}\int \langle \mathcal{L}_1, A_0\rangle\; dx\right)=K_0A_1+C_1,
\end{split}
\end{equation*}
where $C_1$ is a constant of integration. If we set $K_1=C_1A_0^{-1}$, a similar
computation gives:
\begin{equation*}
\begin{split}
B_2 &=-\frac{1}{2\mathcal{L}_0}\int\left( \langle \mathcal{L}_1, B_1\rangle + \mathcal{L}_2, B_0\rangle\right)\; dx\\
 &=-\frac{1}{2\mathcal{L}_0}\int\left( \langle \mathcal{L}_1, K_0A_1+K_1A_0\rangle + \mathcal{L}_2, K_0A_0\rangle\right)\; dx\\
&=K_0\left( -\frac{1}{2\mathcal{L}_0}\int \langle \mathcal{L}_1, A_1\rangle+\langle \mathcal{L}_2, A_0\rangle\; dx\right)+K_1\left( -\frac{1}{2\mathcal{L}_0}\int \langle \mathcal{L}_1, A_0\rangle\; dx\right)\\
&=K_0A_2+K_1A_1+C_2=K_0A_2+K_1A_1+K_2A_0
\end{split}
\end{equation*}
where $C_2$ is the constant of integration and $K_2:=C_2A_0^{-1}$.
We have the following \underline{claim}: there exist constants $K_0$, $K_1$, $\ldots$, $K_s$ such that
\equ\label{induction_Bcoeficients} 
B_j(x) =\sum_{i=0}^j K_{j-i}A_i(x)
\uqe
for $j\leq s$.

By induction, assume the claim is true for $s\leq k$. The recursion formula (\ref{recursion_formula}) for $B_{k+1}$ is
\be\label{Brecursion_formula}
B_{k+1}  =  -\;\frac{1}{2\mathcal{L}_0}\left( \int      \sum\limits_{i = 1}^{d}    \langle  \mathcal{L}_i , B_{k+1-i} \rangle       \; dx  - \frac{1}{2}B_{k+1-d}''\right).
\ee

Hence, using the induction hypothesis of the  \underline{claim} for $s\leq k$,
\equ\label{beelinearsum}
\begin{split}
\sum\limits_{i = 1}^{d}    \langle  \mathcal{L}_i , B_{k+1-i} \rangle &= \sum\limits_{i = 1}^{d}    \langle  \mathcal{L}_i , \sum_{j=0}^{k+1-i} K_{k+1-i-j}A_j \rangle\\
&= \sum_{i=1}^{d}\sum_{j=0}^{k+1-i} K_{k+1-i-j}\langle \mathcal{L}_i, A_j\rangle
\end{split}
\uqe

Descending ordering by the index $k+1-i-j$, we obtain that the right-hand-side of
(\ref{beelinearsum}) can be written as (notice that we are taking the change of variable $l=i+j-1$):

\equ\label{beepartialsum}
K_k\langle \mathcal{L}_1, A_0\rangle+K_{k-1}\left(\langle \mathcal{L}_1, A_1\rangle+\langle \mathcal{L}_2, A_0\rangle\right)+\ldots=\sum_{l=0}^k K_{k-l} \left(\sum_{i=0}^d \langle  \mathcal{L}_i, A_{l+1-i} \rangle \right)
\uqe

But, also
\equ\label{beederivative}
\begin{split}
B^{''}_{k+1-d} &=\left(\sum_{i=0}^{k+1-d} K_{k+1-d-i}A_i\right)^{''} \\
&=\sum_{i=1}^{k+1-d}  K_{k+1-d-i}A_i^{''}=\sum_{l=0}^k K_{k-l} A_{l+1-d}^{''} 
\end{split}
\uqe

Combining expressions (\ref{beelinearsum}-\ref{beederivative}), 
equation (\ref{Brecursion_formula}) becomes
\begin{equation*}
\begin{split}
B_{k+1} &= \sum_{l=0}^k K_{k-l}\left( -\;\frac{1}{2\mathcal{L}_0}\left( \int      \sum\limits_{i = 1}^{d}    \langle  \mathcal{L}_i , A_{l+1-i} \rangle       \; dx  - \frac{1}{2}A_{l+1-d}''\right)\right) \\
&= \left(\sum_{l=0}^k K_{k-l}A_{l+1}\right)+C_{k+1}
\end{split} 
\end{equation*}

But, setting $K_{k+1}:=C_{k+1}A_0^{-1}$, we have:
$$
B_{k+1}=\sum_{l=0}^{k+1} K_{k+1-l}A_l
$$
as needed for proving the \underline{claim}, equation (\ref{induction_Bcoeficients}).

Now, using the \underline{claim}, we can express  $\psi_n$ in terms of $\phi_n$, $\phi_{n-1}$, $\ldots$, $\phi_0$:
\equ
\begin{split}
\psi_n &= \sum_{k=0}^n B_k\lambda^{n-k}
= \sum_{k=0}^n \left( \sum_{i=0}^k K_{k-i} A_i\right) \lambda ^{n-k}\\
&= \sum_{k=0}^n \left( \sum_{j=0}^k K_j A_{k-j}\right) \lambda ^{n-k} \mbox{ (setting  $j=k-i$)}\\
&= \sum_{j=0}^n \sum_{k=j}^n K_j A_{k-j} \lambda ^{n-k} \mbox{ (exchanging the order of summing)}\\
&= \sum_{j=0}^n K_j\left( \sum_{k=j}^n A_{k-j} \lambda ^{n-k}\right)
= \sum_{j=0}^n K_j\left( \sum_{l=0}^{n-j} A_{l} \lambda ^{n-l-j}\right) \mbox{ (setting $l=k-j$)}\\
&= \sum_{j=0}^n K_j \phi_{n-j}
\end{split}
\uqe 

\item[B)]
From the proof of part A), we can see that $\psi_n= \sum_{j=0}^n K_j \phi_{n-j}$
if, and only if, the coefficients $B_j$  of $\psi_n$ relates to the coefficients of
$\phi_n$, by the assertion of the equation (\ref{induction_Bcoeficients}), for $0\leq j\leq n$.

Now, the definition of $\mathcal{B}_{n,s}$ in (\ref{d_conditions}) applied to the coefficients of $\psi_n$ is the equation
$$
\mathcal{B}_{n,s}:= B_{n-s}''' - 2  \sum\limits_{i = d-s}^{d}    \langle  \mathcal{L}_i , B_{d + n-s-i} \rangle  
$$

But, 
\equ\label{sumConditionBns1}
\begin{split}
\sum_{i = d-s}^{d}    \langle  \mathcal{L}_i , B_{d + n-s-i} \rangle   &= \sum_{i = d-s}^{d}    \langle  \mathcal{L}_i , \sum_{j=0}^{d+n-s-i} K_{d+n-s-i-j}A_j \rangle   \\
&=\sum_{i = d-s}^{d}  \sum_{j=0}^{d+n-s-i}  K_{d+n-s-i-j}\langle \mathcal{L}_i, A_j\rangle\\
&=\sum_{l=0}^n K_l\left(\sum_{i = d-s}^d \langle \mathcal{L}_i, A_{n+d-s-i-l}\rangle\right),
\end{split}
\uqe
where the last equation follows by setting $l=d+n-s-i-j$ and ordering by $l$.

Also,
$
B_{n-s}^{'''}=\sum_{i=0}^{n-s}K_{n-s-i}A_i^{'''}=\sum_{l=0}^{n-s}K_{l}A_{n-s-l}^{'''}.
$
Now, if $l>n-s$, then $A_{n-s-l}=0$ (because $n-s-l<0$). Therefore, we can write
\equ\label{sumCondtionBns2}
B_{n-s}^{'''}=\sum_{l=0}^{n}K_{l}A_{n-s-l}^{'''}
\uqe

Thus, combining equations (\ref{sumConditionBns1}) and (\ref{sumCondtionBns2}), we obtain:
\begin{equation*}
\begin{split}
\mathcal{B}_{n,s} &=\sum_{l=0}^{n}K_{l}A_{n-s-l}^{'''}-2\sum_{l=0}^n K_l\left(\sum_{d-s}^d \langle \mathcal{L}_i, A_{n+d-s-i-l}\rangle\right)\\
&=\sum_{l=0}^n K_l\left(A_{n-s-l}^{'''}-2\sum_{d-s}^d \langle \mathcal{L}_i, A_{n+d-s-i-l}\rangle\right)\\
&=\sum_{l=0}^n K_l\mathcal{A}_{n-l,s} 
\end{split}
\end{equation*}
Then, equation (\ref{linear_combination_of_conditions}) holds, since 
%
%
%
%
%
$\mathcal{A}_{n-l,s}=0$ is the $s^{\text{th}}$ condition for the normalized soliton  $\phi_{n-l}$. Hence, $\mathcal{B}_{n,s}=0$ is the $s^{\text{th}}$ condition for the  soliton  $\psi_{n-l}$, as we just proved.

\une 
\frp


\section{Geometry of the $n$-soliton.}\label{sect:hyperelliptic_curve}

In this section, we construct the hyperelliptic curve corresponding to the $n$-soliton solution related to the eigenvalue  problem, equation  (\ref{inicial}),
where  $\hat{\mathcal{L}}$ is  of degree $d$ in  $\lambda$. We obtain a hyperelliptic curve,
$$
\mathcal{C}_n:\frac{1}{2}Y^2  =   -   \mathcal{H}_n(X)
$$
where $\mathcal{H}_n$ is a polynomial of degree $2n+d$ in the variable $X$, with constant coefficients. Thus, the genus of the curve  $\mathcal{C}_n$ is $n + (d-1)/2$ or $n + (d-2)/2$, depending on if $d$ is odd or even, respectively. The factorization of the $n$-soliton $\phi_n=\sum_{i=0}^nA_i\lambda^{n-i}$ will give $n$ different roots $\lambda_k(x)$ depending on $x$, such that the points  $P_k= (X_k, Y_k) = (\lambda_k,\phi_n' (x, \lambda_k))$, for $1\leq k\leq n$, will represent $n$ points on the curve. 

Hence, if we vary $x$, a soliton solution defines $n$ real curves on $\mathcal{C}_n$.
These curves on the tangent space define a system of linear differentials. We will also compute  these differentials in this section. 


\subsection{Computation of the function $\mathcal{H}_n(\lambda)$.}

Define 
\equ\label{definition_H}
%
\mathcal{H}(\phi)
:=\int \phi \left(\phi'''-2\langle \mathcal{\hat{L}}, \phi \rangle \right)dx
=\int \phi\left(\phi'''-4\mathcal{\hat{L}}\phi' -2\mathcal{\hat{L}'}\phi\right)dx
\uqe
Thus, integrating: 
$
\int\phi\phi'''dx=\phi\phi''-\int\phi'\phi''dx=\phi\phi''-\frac{1}{2}\left(\phi'\right)^2
$
and \\
$
-\int (4\phi\mathcal{\hat{L}}\phi' +2\phi\mathcal{\hat{L}'}\phi)dx 
=-2\int(\mathcal{\hat{L}}(2\phi\phi')+\mathcal{\hat{L}'}\phi^2)dx
= -2\int \left(\mathcal{\hat{L}}\phi^2\right)'dx=-2\mathcal{\hat{L}}\phi^2,
$
we obtain
\equ
\mathcal{H}(\phi)  = \phi\phi" - \frac{1}{2}(\phi')^2 - 2 \hat{\mathcal{L}}\phi^2.
\uqe
Observing the leading term, we obtain the following result.

\lmm\label{H_n_higest_term}
If $\phi$ is an $n$-soliton and $deg(\mathcal{\hat{L}})=d$, then
$\mathcal{H}(\phi)$ is a polynomial in $\lambda$, constant with respect to the variable $x$, of degree $2n+d$ with leading coefficient $-2\mathcal{\hat{L}}_0A_0^2$, where $A_0$ is the leading coefficient of $\phi$.
\mml

Now, if $\phi_n$ is the normalized $n$-soliton of degree $n$, set
the function
\be	\label{H_n_function}
\mathcal{H}_n(\lambda)  :=\mathcal{H}(\phi_n)= \phi_n\phi_n" - \frac{1}{2}(\phi_n')^2 - 2 \hat{\mathcal{L}}\phi_n^2.
\ee

Next, we will give a precise formula to compute $\mathcal{H}_n(\lambda)$ as a polynomial in $\lambda$.

Consider the $n$-soliton solution (\ref{squared_eigenfuntion_n}) and the multiplicative operator (\ref{L_degree_d}), 
where $A_0$ and $\mathcal{L}_0$ are non-zero constants ({\em i.e.}, $x$-independent);  $A_i$ and $\mathcal{L}_j$  are functions of $x$; and $\lambda$ is a constant. 

There are two ways to explicitly compute $\mathcal{H}_n(\lambda)$: 
\enu
\item \textbf{Hard way}: Use directly the expression  for $\mathcal{H}_n(\lambda)  = \phi_n\phi_n" - \frac{1}{2}(\phi_n')^2 - 2 \mathcal{L}\phi_n^2,$. We did it this way with the KdV equation. Details can be found in the proof of \textbf{Theorem 7.6} in Appendix A of \cite{espinola_portillo}. 
\item \textbf{Easy way}: Use bilinearity in the equation
$$
\mathcal{H}_n(\lambda)  =   \int \phi_n\left(  \phi_n'''-2\langle \mathcal{\hat{L}}, \phi_n \rangle   \right) dx + \text{constant of integration.}
$$
Notice that  when solving $ \phi_n'''-2\langle \mathcal{\hat{L}}, \phi_n \rangle=0 $, we obtain:  
\enu
\item the recursive formul\ae  \; to compute the coefficients  $A_i$ of $\phi_n$;
\item additional $d$ conditions of solvability.   
\une
\une

Thus, from the discussion in section \ref{sec:recursion_formulas_conditions}, we have 
\equ\label{treintaytantos}
\phi_n'''-2\langle \mathcal{\hat{L}}, \phi_n \rangle=\sum\limits_{i = 1}^{d} \mathcal{A}_{n,d-i} \lambda^{d-i}, 
\uqe
where $\mathcal{A}_{n,d-i}=0$ are the $d$ conditions in equation  (\ref{d_conditions}).
From the previoues equation and the fact that $\phi_n$ is a polynomial of degree $n$,
\bas
\mathcal{H}_n'(\lambda)  & = &  \phi_n  \left(\phi_n''' - 2 \langle  \mathcal{L}  ,   \phi_n\rangle   \right) \\
          & = &  \sum_{k = 1}^{d+n}  \left(\sum_{\scriptsize\begin{matrix}  i + j= k\\ d\geq i \geq 1\\ n\geq j \geq 0\end{matrix}}   \mathcal{A}_{n,d-i} A_j   \right)  \lambda^{d+n-k}
\eas





Integrating (with respect to $x$) and observing that the highest degree in $\mathcal{H}_n(\lambda) $  is $-2\mathcal{L}_0A_0^2 \lambda^{2n+d}$ (which turns to be the constant of integration by lemma (\ref{H_n_higest_term})), we obtain the following theorem.

\thm\label{theorem_formula_H_n} The function $\mathcal{H}_n(\lambda)$ is a polynomial of degree $2n+d$ with constant coefficients ({\em i.e.}, they are $x$-independent) and is given by the formula
\equ\label{H_n_formula}
\mathcal{H}_n(\lambda) = -2\mathcal{L}_0A_0^2 \lambda^{2n+d}  + \sum_{k = 1}^{d+n}  \left(\sum_{\scriptsize\begin{matrix}  i + j= k\\ d\geq i \geq 1\\ n\geq j \geq 0\end{matrix}}  \int  \mathcal{A}_{n,d-i} A_j  \; dx \right)  \lambda^{d+n-k} .
\uqe

\mht

\textbf{Remark.} Please note that there is a gap ({\em i.e.}, there is no terms) from $\lambda^{n+d}$ and up to $\lambda^{2n+d-1}$.


\subsection{The hyperelliptic curve associated to the solutions.}

Consider the normalized $n$-soliton $\phi_n(x)=\sum_{i=0}^n A_i(x) \lambda^{n-i}$.
We can factorize $\phi_n(x)$ over $\C$ for fixed
values of $x$  as
$$
\phi_n(x;\lambda)=A_0\prod_{i=1}^{n}\left[\lambda-\lambda_i(x)\right].
$$
Here, $\lambda_i(x)$ are the roots of $\phi_n(x)$ which depend on $x$.

Taking the derivative  with respect to $x$ of the previous expression, we obtain
$$
\phi'_n(x;\lambda)=-A_0\sum_{j=1}^n\lambda'_j\prod_{i\neq j}\left[\lambda-\lambda_i(x)\right],
$$
evaluating at $\lambda=\lambda_k(x)$, we finally get
\equ\label{eqg1}
\phi'_n (x, \lambda_k(x))  =  -A_0\lambda'_k\prod_{i\neq k}\left[\lambda_k(x)-\lambda_i(x)\right].
\uqe

Now,  evaluating $\mathcal{H}_n(\lambda)$ in equation (\ref{H_n_function})
at $\lambda=\lambda_k(x)$, and using that $\phi_n\mid_{\lambda=\lambda_k(x)}=0$, we obtain
\equ\label{eqg2}
\mathcal{H}_n(\lambda_k)=-\frac{1}{2}(\phi'_n (x, \lambda_k(x)))^2
\uqe

Now, since $\mathcal{H}_n(\lambda)$ is a polynomial of degree $2n+d$ (Theorem \ref{theorem_formula_H_n}) with constant coefficients with respect to the variable $x$, the equation
$$
\mathcal{H}_n(X)=-\frac{1}{2}Y^2
$$
is the equation of an hyperelliptic curve $\mathcal{H}_n$ of genus $n+(d-1)/2$ if $d$ is odd, or $n+(d-2)/2$ if $d$ is  even; each $P_k= (X_k, Y_k) = (\lambda_k(x),\phi_n'(x,\lambda_k(x))$,  
with $1\leq k\leq n$, represents a point on the curve. Hence, if we vary $x$, a soliton solution defines $n$ real curves on $\mathcal{C}_n$.
The curves on the tangent space define a system of linear differentials. 
Next, we compute these differentials. 

In fact, combining equations (\ref{eqg1}) and (\ref{eqg2}), we have that 
$$
\sqrt{-2\mathcal{H}_n(\lambda_k)} =\phi'_n (x, \lambda_k(x)) = 
 -A_0\lambda'_k\prod_{i\neq k}\left[\lambda_k(x)-\lambda_i(x)\right].
$$
Hence, setting $\mathcal{R}_n=-2\mathcal{H}_n,$ we obtain
$$
\frac{\lambda_k'}{\sqrt{\mathcal{R}_n(\lambda_k)}}=\frac{-A_0^{-1}}{\prod_{i\neq k}\left[\lambda_k(x)-\lambda_i(x)\right]}
$$
Now, the differentials 
$$
\omega_\mu=\frac{X^{\mu-1}dX}{\sqrt{\mathcal{R}_n(X)}}\mbox{ , $1\leq\mu\leq n$}
$$
form a basis of the  space of differentials of the curve $\mathcal{H}_n$.

Evaluating those differentials at the points $P_k$, we obtain
$$
\omega_\mu(P_k)=\omega_\mu(\lambda_k) = \frac{\lambda_k^{\mu-1}\lambda_k'dx}{\sqrt{\mathcal{R}_n(\lambda_k)}}
=\frac{-A_0^{-1}\lambda_k^{\mu-1}}{\prod_{i\neq k}\left[\lambda_k(x)-\lambda_i(x)\right]} 
$$
Adding up over all points $P_k$, and using the main Proposition 
in Appendix B 
in \cite{espinola_portillo}, we get
\equ\label{maindiff}
\sum_{k=1}^{n}\omega_\mu(P_k)=
\begin{cases}
0 & \mbox{ if }1\leq \mu< n \cr
-A_0^{-1} & \mbox{ if }\mu=n.
\end{cases}
\uqe
Other proofs of this fact are found in \cite{newell, miller}. 
These differentials play a fundamental role on integrating the PDEs under consideration. 



\section{Examples}\label{examples}\label{sect:examples}

\subsection{ The Korteweg-deVries (KdV) equation.}\label{subsect:kdv}

For the the Korteweg-deVries equation (\ref{kdv_equation}), we have the associated linear eigenvalue problem (\ref{inicial}), 
where  $\hat{\mathcal{L}}=\lambda-q$ is the  Schr\"odinger operator. Equation (\ref{rec:4}) becomes the 
linear differential equation: 
\equ\label{eq1}
\phi'''+4q\phi'+2q'\phi=4\lambda \phi',
\uqe
 If we set the
linear differential operator
\equ\label{kdv-operator}
\mathcal{M}:=\frac{d^3}{dx^3}+4q\frac{d}{dx}+2q',
\uqe
 the differential equation (\ref{eq1}) now reads
\equ\label{kdv-phi}
\mathcal{M}(\phi)=4\lambda\phi'.
\uqe
Set
\equ\label{equ2}
\phi_n(x,t)=F_{-1}(4\lambda)^n+F_{0}(4\lambda)^{n-1}+\cdots+F_{n-1}
\uqe
where $F_{-1}$ is constant, and $F_i(x)$ for $-1\leq i <n-1$ are {\em implicit} functions of $x$. (The factor "4" just helps to normalize the conserved quantities). Notice the  shift in the subscript as opposed to equation (\ref{squared_eigenfuntion_n}).
Now, (\ref{kdv-phi}) implies the recursion formula:
\equ\label{kdv:recursion}
F_{j}=\int \mathcal{M}(F_{j-1}) dx, \makebox[1.5in]{ for $j=0$,$1$,$\ldots$,$n-1$.}
\uqe
and the condition:
\equ\label{condition_F_constant}
F_{n}:=\int \mathcal{M}(F_{n-1}) dx \quad  \makebox[1.0in]{ is a constant.}
\uqe
The function $F_n(x)$ is  the {\em $n^{\text{th}}$conservation density} of the KdV equation.

If we normalize setting $F_{-1}=\frac{1}{2}$, the fist few coefficients of $\phi_n$  are: 
\begin{eqnarray}
F_{0}  & = & q\mbox{;     } \nonumber  \\
 F_{1} & =  &   q''+3q^2\mbox{;     }\\ \label{KdV_conserved_densities}
 F_{2}  & = &  q^{(4)}+10qq''+5(q')^2+10q^3\mbox{,} \nonumber
\end{eqnarray}
which turn to be the conserved densities of the KdV equation. 

Using the results of section \ref{sect:recursion}, we can recover the formul\ae  \; just presented. The recursion formula, equation (\ref{recursion_formula}) is 
$$
A_{k}'  =  -\;\frac{1}{2\mathcal{L}_0}\left(      \sum\limits_{i = 1}^{d}    \langle  \mathcal{L}_i , A_{k-i} \rangle         - \frac{1}{2}A_{k-d}'''\right),
$$
and the extra condition, equation (\ref{d_conditions_alternative}), is $\mathcal{A}_{n,s}= 0$, with 
$$
\mathcal{A}_{n,s}= 4\mathcal{L}_0A'_{n+d-s}+2\sum\limits_{i = 1}^{d-s-1} \langle  \mathcal{L}_i , A_{d + n-s-i} \rangle .
$$

For the KdV equation, the linear operator is $\hat{\mathcal{L}} = \lambda - q(x)$, and the $n$-soliton solution 
is $\phi_n(x) = A_0 \lambda^n + A_1(x)\lambda^{n-1} + \dots A_n(x)$. Therefore, $d= 1, \mathcal{L}_0 = 1, \mathcal{L}_1(x)  =  - q(x)$. Then, the recursion formula becomes:
$$
A_{k}'  =  \frac{1}{2}    q'(x)  A_{k-1}  +   q(x) A_{k-1}'       + \frac{1}{4}A_{k-1}'''.
$$
and we recover the recursion formula in  equation (\ref{kdv:recursion}).

For the extra-conditions, we just have one, with $s=0$ ($0\leq s \leq d-1 = 0$), 
{\em i.e.},
$$
 4\mathcal{L}_0A'_{n+1}+2\sum\limits_{i = 1}^{0} \langle  \mathcal{L}_i , A_{1 + n - i} \rangle = 0.
$$
Therefore,  the unique condition is that $A_{n+1} = 4 F_{n}$ is constant. Same as condition (\ref{condition_F_constant}). 










\dfn \label{defpsi}
A function of the form
$$
\psi_n=A_{-1}(4\lambda)^n+A_0(4\lambda)^{n-1}+\cdots+A_{n-2}(4\lambda)+A_{n-1},
$$
(with $A_{-1}$ constant and $A_i(x)$ are functions of $x$ for $0\leq i \leq n-1$)
solution of $B(\psi_n)=4\lambda\psi_n'$, is called a {\em KdV $n$-soliton.}
The {\em normalized KdV $n$-soliton}, denoted by $\phi_n$, is the KdV $n$-soliton  
setting: $A_k=4^{n-k}F_{k-1}$. 
\nfd

Theorem \ref{all_N_solitons} in the KdV case specializes as follows.

\thm\label{all_KdV_N_solitons}
Each KdV $n$-soliton $\psi_n$ can be written as a linear combination of the normalized solitons: $\phi_n$, $\phi_{n-1}$, $\ldots$, $\phi_0$.
Moreover, a linear combination
$$\psi_n=\alpha_n\phi_n+\alpha_{n-1}\phi_{n-1}+\ldots+\alpha_0\phi_0$$ 
(with $\alpha_i$ constant and $\alpha_n\neq 0$) is a KdV $n$-soliton if, and only if,
$$\alpha_nF_n+\alpha_{n-1}F_{n-1}+\cdots+\alpha_0F_0$$ 
is constant.

\mht

Also, Theorem \ref{theorem_formula_H_n} has its KdV version.

\thm \label{thm2}
The function $\mathcal{H}_n(\lambda):=\int\phi_n B(\phi_n)dx -2\lambda\phi_n^2$ is given by the formula:
\equ
\begin{split}
\mathcal{H}_n(\lambda) &=-\frac{(4\lambda)^{2n+1}F_{-1}^2}{2}+(4\lambda)^{n}F_{-1}F_n+(4\lambda)^{n-1}\left[F_0F_n-\int F_0'F_ndx \right]+\\
&\makebox[5mm]{ }+(4\lambda)^{n-2}\left[F_1F_n-\int F_1'F_ndx \right]+\cdots+\left[F_{n-1}F_n-\int F_{n-1}'F_n dx \right].
\end{split}
\uqe
\mht

\prf
Since $d=1$ and $\mathcal{L}_0=1$ for the KdV equation, Theorem \ref{theorem_formula_H_n} implies that
\equ\label{Hn_kdv_1}
\mathcal{H}_n(\lambda)=-2A^2_0\lambda^{2n+1}+ \sum_{k = 1}^{n+1}  \left(  \int  \mathcal{A}_{n,0} A_{k-1}  \; dx \right)  \lambda^{n+1-k}.
\uqe
But, in notation of section \ref{sect:hyperelliptic_curve}, we assume that
$\phi_n=\sum_{i=0}^nA_i\lambda^{n-i}$, while for the KdV example in this section, we use
the classical representation of $\phi_n$ as the sum 
$$\phi_n=F_{-1}(4\lambda)^n+F_0(4\lambda)^{n-1}+\cdots+F_{n-2}(4\lambda)+F_{n-1}.$$
Equating coefficients, we have the relation among them: $A_j=4^{n-j}F_{j-1}$.
Hence, the leading term of the polynomial $\mathcal{H}_n(\lambda)$ is
$$-2A^2_0\lambda^{2n+1}=-2\left(4^nF_{-1}\right)^2\lambda^{2n+1}=-\frac{(4\lambda)^{2n+1}F_{-1}^2}{2}.$$
Now, using that $\mathcal{A}_{n,0}=4A'_{n+1}$ (by the alternative form in equation (\ref{d_conditions_alternative})), and that $F_n'=4A'_{n+1}$ and $A_{k-1}=4^{n+1-k}F_{k-2}$, we obtain that the sum in (\ref{Hn_kdv_1}) becomes
$$\sum_{k = 1}^{n+1}  \left(  \int  F_{n}' 4^{n+1-k}F_{k-2}  \; dx \right)  \lambda^{n+1-k}=\sum_{j = 0}^{n}  \left(  \int  F_{n}' F_{j-1}  \; dx \right)  (4\lambda)^{n-j}.
$$
Now, the result follows after integration by parts.

\frp

\subsubsection{Examples.} Here, we show the normalized $n$-soliton solutions for $n= 0, 1, 2$ and $3$:
\bas
\phi_0(\lambda) & = & \frac{1}{2} \\
\phi_1(\lambda) & = &  2\lambda  + q\\
\phi_2(\lambda) & = &  8 \lambda^2   +  4q\lambda  + q" + 3q^2\\
\phi_3(\lambda) & = & 32\lambda^3+16q\lambda^2+4(q" + 3q^2)\lambda+q^{(4)}+10qq''+5(q')^2+10q^3  
\eas

And, for $n= 0, 1$ and $2$, their corresponding
hyperelliptic curves are
\bas
\mathcal{H}_0(\lambda) & = & -\frac{1}{2}(\lambda - q)  \\
\mathcal{H}_1(\lambda) & = &  -8\lambda^3  + 2\lambda(q" + 3q^2) +   \left[  q(q" + 3q^2) - \left( \frac{1}{2}(q')^2 +  q^3\right)   \right]\\
\mathcal{H}_2(\lambda) & = &  -128 \lambda^5 +8\lambda^2 F_2 + 4\lambda K_2   +  L_2 ,
\eas
where
\bas
F_2 & = & q^{(4)}+10qq''+5(q')^2+10q^3,    \\
K_2 & = &  qq^{(4)}-q'q^{(3)}+10q^2q^{(2)}+\frac{1}{2}(q^{(2)})^2+10q(q')^2+\frac{25}{2}q^4,\\
L_2 & = & q^{(4)}(q^{(2)}+3q^2)-q^{(3)}(\frac{1}{2}q^{(3)}+6qq')+q^{(2)}(8q^{(2)}q+6(q')^2+30q^3)+18q^5. 
\eas


\subsection{The Nonlinear Schr\"odinger (NLS) equation.}\label{subsect:nls}

The stationary version of the NLS equation (\ref{nls_equation}) is
$$
  \frac{1}{2}\frac{\partial^2 q}{\partial x^2}  + \sigma |q|^2 q  =   \omega q,
$$
where $\sigma = \pm 1$ is the focusing/defocusing parameter and $\omega$ is a constant. 
By the work of Kamchatnov, Kraenkel and Umarov \cite{Kamchatnov_Kraenkel_Umarov_1, Kamchatnov_Kraenkel, Kamchatnov_Kraenkel_Umarov_2}, we know that the multiplicative operator $\hat{\mathcal{L}}$ of the NLS equation 
is a polynomial of degree $d=2$ in $\lambda$:
\equ
\begin{split}\label{nls:equation:E:F}
\hat{\mathcal{L}} &=-\left(\lambda-\frac{iq_x}{2q}\right)^2-\sigma\|q\|^2-\left(\frac{q_x}{2q}\right)_x\\
&= -\lambda^2+E\lambda+F
\end{split},
\uqe
with 
\equ\label{Eexpresion}
E=\frac{iq_x}{q} 
\uqe
and 
\equ\label{Fexpresion}
F=-\frac{1}{4}E^2-\sigma\|q\|^2+\frac{i}{2}E'. 
\uqe

Hence, we have the following theorem. 
\thm\label{NLS_theorem}
The following statements are true.
\enu
\item[(i)]  $\mathcal{H}_n(\lambda)$ is a polynomial in $\lambda$ of degree $2n + d = 2n+2$. 
\item[(ii)] The recursion formula is
\equ \label{recursionls:5} 
A_{j}=\frac{1}{2}\int\left[\langle E, A_{j-1}\rangle+\langle F, A_{j-2}\rangle\right] dx-\frac{1}{4}A_{j-2}''
\uqe
for $j=1,\ldots,n$ (assuming $A_{-1}=0$).

\item[(iii)] The following two conditions hold:
\enu
\item[] \textbf{Condition A}, corresponding to $s=1$:
\equ\label{condition:A:differential}
A_{n+1}'=0
\uqe

\item[]\textbf{Condition B}, corresponding to $s=0$:
\equ\label{condition:B:differential}
-4A'_{n+2}+2A_{n+1}E'=0
\uqe 

\une
\une
\mht

\prf
The conditions stated in the Theorem are the alternative conditions $\mathcal{B}_{n,s}=0$ in equation (\ref{d_conditions_alternative}) for the NLS equation. In this case, we have $d=2$, $\mathcal{L}_0=-1$, $\mathcal{L}_1=E$, and thus, the conditions are:
$$\mathcal{B}_{n,1}=-4A'_{n+1}=0$$
and
$$\mathcal{B}_{n,0}=-4A'_{n+2}+2\langle E,A_{n+1} \rangle=-4A'_{n+2}+2A_{n+1}E'=0;$$
hence $A_{n+1}$ is constant, by condition $\mathcal{B}_{n,1}=0$.
\frp




Theorem \ref{all_N_solitons} in the NLS example can be stated as follows.

\thm\label{all_NLS_N_solitons}

Each NLS $n$-soliton solution $\psi_n$ can be written as a linear combination of the normalized solitons: $\phi_n$, $\phi_{n-1}$, $\ldots$, $\phi_0$.
Moreover, consider a linear combination
$$\psi_n=\alpha_n\phi_n+\alpha_{n-1}\phi_{n-1}+\ldots+\alpha_0\phi_0,$$ 
where $\phi_i$ is the normalized NLS $i$-th soliton, $\alpha_i$  is constant, and $\alpha_n\neq 0$; then $\psi_n$ is a NLS $n$-soliton if, and only if,
\enu
\item 
$\alpha_nA_{n+1}+\alpha_{n-1}A_{n}+\cdots+\alpha_0A_1$ is constant.
\item 
$\sum_{j=0}^n \alpha_j(-4A_{j+2}+2A_{j+1}E)$ is constant.
\une
\mht

\prf
The result follows after integration of conditions in theorem \ref{NLS_theorem}.
\frp




 {\theorem The function $\mathcal{H}_n(\lambda)$ for the NLS equation.} 
For the NLS equation, the function $\mathcal{H}_n(\lambda)$, which defines the hyperelliptic curve for the $n$-soliton solution, is given by the following formula:
$$
\mathcal{H}_n(\lambda) = 8\lambda^{2n+2}  - 8 A_{n+1}\lambda^{n+1} + \sum\limits_{i=1}^n \left\{ \int  \left[  - 4A_i A_{n+1}' + A_{i-1}B_n  \right]dx\;  \lambda^{n+1-i} \right\} + \int A_nB_n\;dx, 
$$
where $B_n = A_n''' - 2\langle \mathcal{L}_2, A_n\rangle$.

\begin{proof}
Extending the recursive formula (\ref{recursionls:5}) to $j= n+1$,  Condition A becomes:
\begin{eqnarray*}
0 = A_{n-1}''' - 2\left[  \langle E, A_{n}\rangle +  \langle F, A_{n-1}\rangle \right] &=& -4 \left\{ \frac{1}{2}  \left[  \langle E, A_{n}\rangle +   \langle F, A_{n-1}\rangle \right] -   \frac{1}{4} A_{n-1}'''   \right\}  \\
      & = &     -4 A_{n+1}',
\end{eqnarray*}
{\em i.e.}, $A_{n+1}$ is constant. 
In this case $A_{n+1}$ is the new coefficient for the ($n+1$)-soliton $\phi_{n+1}(x;\lambda)$. It also represents the $(n+1)^{\text{th}}$ conservation density of the NLS equation. 

Hence, by equation (\ref{treintaytantos}), 
$$
\phi_n'''-2\langle \mathcal{\hat{L}}, \phi_n \rangle=  - 4 A'_{n+1}\lambda  + B_n.
$$
and
\bas
\mathcal{H}_n'(\lambda) & = &  \phi_n(\phi_n'''-2\langle \mathcal{\hat{L}}, \phi_n \rangle) = \left( 2\lambda^n  + \sum\limits_{i= 1}^n A_i\lambda^{i-1}\right) \Big(  - 4 A'_{n+1}\lambda  + B_n\Big) \\
& = & - 8A_{n+1}' \lambda^{n+1} +  \sum\limits_{i= 1}^n (-aA_iA_{n+1}' + A_{i-1}B_n)\lambda^{n+1-i}  + A_nB_n.
\eas

Notice that the leading coefficient of $\phi_n$ is $-2\mathcal{L}_0=2$, by the definition of the normalized solitons.
We obtain the result integrating. We just have to notice that the highest degree term in $ \mathcal{H}_n(\lambda) $ is $8\lambda^{2n+2}$, which follows from the highest term in $\lambda$ in the summand $-2\hat{\mathcal{L}}\phi_n^2$ in equation (\ref{H_n_function}). This term is the constant of integration. 

\end{proof}


\subsubsection{Computations of some normalized $n$-soliton solutions}

As in the KdV case, we compute some normalized $n$-soliton solutions for the NLS equation, using
the recursion formula (\ref{recursionls:5}).  

The normalized $0$-soliton is 
$$
\phi_0:=A_0=2.
$$
Thus, 
$
A_1 =\frac{\int \langle E, 2\rangle dx}{2} 
 = \frac{1}{2}\int 2 E_x dx=E+C
$
with $C$ a constant.
If we set $C=0$, we obtain the normalized $1$-soliton
$$\phi_1:=2\lambda+E,$$
with 
$$
\displaystyle A_1 = E= i\frac{q_x}{q}
$$
being a constant for the 1-soliton solution. 

Continuing with the recursion, we obtain
$
A_2 =\frac{1}{2}\int\left[\langle E, E\rangle+\langle F, 2\rangle\right]
$
\\
$
=\frac{1}{2}\int\left[3EE'+2F'\right]
=\frac{3}{4}E^2+F+C. 
$
%
Taking again $C=0$, we get
$$\phi_2:=2\lambda^2+E\lambda+ \left( \frac{3}{4}E^2+F  \right).$$
We can check that  
$$
\displaystyle A_ 2 = \frac{3}{4}E^2+F   = -\frac{1}{2}  \frac{q_{xx}}{q}  - \sigma|q|^2
$$
 is a constant for the 2-soliton solution. See section \ref{sect:1_soliton_NLS}.

Using the recursion formula (\ref{recursionls:5}), and considering  all constants of integration equal to zero, we compute the normalized $3$-soliton solution:
$$
\phi_3:=2\lambda^3+E\lambda^2+ \left(\frac{3}{4}E^2+F\right)\lambda + \left( \frac{5}{8}E^3+\frac{3}{2}FE-\frac{1}{4}E''\right)
$$
with  
$$
\displaystyle A_3 = \frac{q_{xxx}}{q}  +  6\sigma|q|^2\frac{q_x}{q}
$$
constant for the 3-soliton solution. See section \ref{sect:2_soliton_NLS}. 

Similarly, the $4$-soliton solution is:
\begin{align*}
\phi_4:= & 2\lambda^4+E\lambda^3+\left(\frac{3}{4}E^2+F\right)\lambda^2+\left(\frac{5}{8}E^3+\frac{3}{2}FE-\frac{1}{4}E''\right)\lambda \cr
& +\left(\frac{35}{64}E^4+\frac{15}{8}E^2F+\frac{3}{4}F^2-\frac{5}{16}\left(E'\right)^2-\frac{5}{8}EE''-\frac{1}{4}F''\right) .
\end{align*}




\subsubsection{Succesive derivatives of $E=\frac{iq_x}{q}$.}\label{Sec:E:derivatives}

The first derivative of $E$ is
$
E'=\frac{iq_{xx}}{q}-i\frac{q_x^2}{q^2}.
$
Setting 
\equ\label{E2:definition}
E_{(2)}=\frac{iq_{xx}}{q},
\uqe
 we can write the first derivative as:
\equ\label{E:first:derivative}
E'=E_{(2)}+iE^2
\uqe
Now, define:
\begin{equation}\label{E_(n)}
E_{(n)}:=\frac{iq^{(n)}}{q}.
\end{equation}
Hence, we easily compute:
$
E_{(n)}'=E_{(n+1)}+iE_{(n)}E. 
$
%
Using this notation, we can easily compute higher order derivatives of $E$.
For example,
$
E'' =E_{(2)}'+2iE E' =E_{(3)}+iE_{(2)}E+2iE\left(E_{(2)}+iE^2\right),
$
{\em i.e.},
%
\equ\label{E:second:derivative}
E'' =E_{(3)}+3iE_{(2)}E-2E^3
\uqe
and
%
$
E''' =E_{(3)}'+3i\left(E_{(2)}E\right)'-6E^2E'
=E_{(4)}+4iE_{(3)}E-12E^2E_{(2)}-6iE^4+3iE_{(2)}^2.
$

\subsubsection{0-soliton solution for the NLS equation}

We have that $\phi=A_0$ is a constant. Since $\hat{\mathcal{L}} =-\lambda^2+E\lambda+F$ thus, Condition A (equation (\ref{condition:A:differential})) gives
$
\langle E, A_0\rangle=A_0E'=0. 
$
Hence,   
$E= iq_x/q=i\left(\ln{q}\right)'=k$ is constant.
Then, $q$ satisfies the linear equation
$$q_x=-ikq,$$
which is the $0^{\text{th}}$  equation in the NLS Lax hierarchy. 
Thus,  it follows that 
$q=Ce^{-ikx},$ with $C$ and $k$ constants.

Condition B in (\ref{condition:B:differential}) is 
$
\langle F, A_0\rangle=A_0F'=0.
$
Hence, $F$ is constant and 
\begin{equation*}
F =-\frac{1}{4}E^2-\sigma\|q\|^2+\frac{i}{2}E'  
=-\frac{1}{4}k^2-\sigma\|q\|^2.
\end{equation*}
Hence, $\|q\|^2$ is constant, because $F$, $E$ and $\sigma$ are constants.
This is the first conserved density of the NLS equation. 
Thus, we can conclude that $q=Ce^{-ikx}$ with $k\in\R$.

\subsubsection{1-soliton solution for NLS equation}\label{sect:1_soliton_NLS}



Condition A in (\ref{condition:A:differential}) implies  that
$A_2=\frac{3}{4}E^2+F$ is constant. 
Now, using the expresion of $F$ in (\ref{Fexpresion}), we conclude that
$$
A_2 = \frac{1}{2}E^2-\sigma \|q\|^2+\frac{i}{2}E' = \text{constant.}
$$

But, from (\ref{E2:definition}) and (\ref{E:first:derivative}) and from defining
the constant $\omega=-A_2$,  we obtain the stationary nonlinear Schr\"odinger equation:
\equ\label{Stacionary-NLS}
\sigma\|q\|^2q+\frac{1}{2}q_{xx}=\omega q .
\uqe








Condition B in (\ref{condition:B:differential}) implies that  $-4A_3+2A_2E = \Omega$ is constant. After substitution of the values of $A_2$ and $A_3$, it becomes:
$-E^3-4FE+E'' = \Omega.$
Now using (\ref{Fexpresion}) and (\ref{E:second:derivative}), we simplify 
to:
\equ\label{Second:Condition:1-soliton:simplified}
4\sigma \|q\|^2q_xq+q_{xxx}q-q_{xx}q_x=\Omega q^2.
\uqe
%
From (\ref{Stacionary-NLS}), we obtain:
$
q_{xx}=2\omega q-2\sigma \bar{q}q^2,
$
and taking its derivative,
$
q_{xxx}=2\omega q_x-2\sigma \bar{q}_x q^2-4\sigma\|q\|^2q_x.
$
Substituting these expressions for $q_{xx}$ and $q_{xxx}$ to reduce the order of the derivatives  in (\ref{Second:Condition:1-soliton:simplified}), we finally obtain that
\equ\label{Conservation:law:1-soliton}
\bar{q}q_x-\bar{q}_xq=\frac{\sigma\Omega}{2}
\uqe
where  $\Omega$ a constant. This is the second conserved density for the NLS equation.


\subsubsection{2-soliton for NLS}\label{sect:2_soliton_NLS}






Condition A for the 2-soliton becomes $A_3=\frac{5}{8}E^3+\frac{3}{2}FE-\frac{1}{4}E''$ equals  a constant.
Multiplying by $-4$ and substituting the values of $E''$, $E'$ and $F$, and using
(\ref{E:first:derivative}), (\ref{E:second:derivative}) and (\ref{Fexpresion}), we obtain 
$6\sigma\|q\|E+E_{(3)}=i\omega_2,$ 
where $\omega_2$ is constant.
In terms of $q$, we have the following condition:
\equ\label{condition:A2:q}
q_{xxx}+6\sigma\|q\|^2q_x=\omega_2q,
\uqe
which is is the complex modified Korteweg-deVries (mKdV) equation, which is the second  equation in the NLS Lax hierarchy.

Now,  Condition B for the 2-soliton is $-4A_4+2A_3E = $ constant, which in terms of $E$, $F$ and their derivatives, is
$$
\frac{-15}{16}E^4-\frac{9}{2}FE^2-3F^2+\frac{5}{4}(E')^2+2EE''+F''  = \text{constant.}
$$
Using (\ref{Fexpresion}) to substitute $F$ and $F''$ and express only in terms of $E$ and its derivatives, we obtain that the following expression is a constant:
$$
\frac{3}{2}(E')^2+\frac{3}{2}-\frac{3i}{2}E^2+\frac{i}{2}E'''+3\sigma\|q\|^2E^2+3i\|q\|^2E'-3\|q\|^4-\sigma\|q\|^2_{xx}.
$$
From the expressions of $E'$, $E''$ and $E'''$ in subsection \ref{Sec:E:derivatives}, we simplify this second condition to:
$$ 
-\frac{1}{2}EE_{(3)}+\frac{i}{2}E_{(4)}+3i\sigma\|q\|^2E_{(2)}-3\|q\|^4-\sigma\|q\|^2_{xx}=\Omega_2,
$$ 
where $\Omega_2$ is a constant.
In terms of $q$ and its derivatives, we get
\equ\label{condition:B2:q}
\frac{1}{2}q_xq_{xxx}-\frac{1}{2}q_{xxxx}q-3\sigma\|q\|^2 q_{xx}q-3\|q\|^4q^2-\sigma\|q\|^2_{xx}q^2=\Omega_2q^2
\uqe
Using condition (\ref{condition:A2:q}) and its derivative (in order to reduce the order of the derivatives in (\ref{condition:B2:q})), we finally obtain, 
\equ\label{Conservation:law:2-soliton}
\|q_x\|^2-\sigma \|q\|^4-\frac{1}{3}\|q\|^2_{xx}=\frac{\sigma\Omega_2}{3},
\uqe
which is the third conserved density of the NLS equation.


\section{Conclusions}\label{examples}\label{sect:conclusions}

In this paper, we have considered the spectral problem associated to completely  integrable partial differential equations (PDEs) in the sense of Lax pairs theory. The spectral operator is 
assumed to be scalar, linear, multiplicative, 
and of polynomial form of degree $d$ in the spectral parameter $\lambda$. We assume that the solution to the PDE, $q(x)$, is stationary and periodic in $x$. 
We translated the spectral problem into a linear third order differential equation for the associated squared eigenfunctions (this is a standard linearization of the problem). We called $n$-solitons to the polynomial solutions of degree $n$ in $\lambda$ of the squared eigenfunction equations.

Then, we rewrite  the linearized problem  by introducing a bilinear form. Hence, using linear algebra and matching coefficients of same degree, we obtain recursion formul\ae \, to compute the basic $n$-solitons, which we proved generate the solutions to the considered differential equation. We called this formula the recursion formula of the $n$-solitons.

We also discovered $d$ extra conditions which are necessary to be satisfied in order to have solutions to the considered equation. This $d$ conditions give also important  information of the system.
We also show that the $n$-solitons solutions can be parametrized by points on a hyperelliptic curve with genus $n+(d-1)/2$ or $n+(d-2)/2$, depending on if $d$ is odd or even, respectively. We find a formula for the equation of the hyperelliptic curve in terms of the coefficients of the $n$-solitons and the $d$ conditions of the system, which define constants of motion of the system.

As examples, we consider two classical equations, the KdV and the NLS equations, which have multiplicative scalar operators of degre $d=1$ and $d=2$, respectively. 
We found (the recursion formula for) the coefficients  of the $n$-soliton solutions (the squared eigenfuntions), the hyperelliptic curve and  $d$ extra conditions for each case. 
In the KdV case, the extra condition ($d=1$) states  that  the coefficients of the $n$-soliton solutions are constants and they turn to be the conserved densities of the KdV equation. 
For the NLS case, we have two extra conditions ($d=2$). One of them provides us again  the coefficients of the $n$-soliton which are  the conserved densities of NLS. The second condition
represents the corresponding equations of the Lax hierarchy of the NLS equation. 

Finally, these approach can be used in other integrable PDEs, if the pair Lax associated to the PDE is a scalar problem. 







\section*{Aknowledgment.} The authors are very thankful with E. Mart\'inez-Ojeda (UACM) for the comments, suggestions  and useful discussions regarding 
 the present work. 
 
 Esp\'{\i}nola-Rocha dedicates this work to the memory of his teacher and mentor, Prof. A.A. Minzoni, who passed away on July 1st, 2017. 
 The first beauties of integrable systems that he learned, were taught by Minzoni.
 
Portillo-Bobadilla dedicates this work to his wife and son, his supportive family. And also to the beautiful city of Guanajuato!


\end{document}